\documentclass[12pt, oneside, a4paper]{article}

\usepackage{amsmath}
\usepackage{amsfonts}
\usepackage{amssymb}
\usepackage{amsthm,mathrsfs}
\usepackage{enumerate}
\usepackage{graphicx}
\newtheorem{theorem}{Theorem}[section]
\newtheorem{lemma}[theorem]{Lemma}

\theoremstyle{definition}

\title{\textbf{The character graph of a finite group is perfect}}
\author{Mahdi Ebrahimi\footnote{ m.ebrahimi.math@ipm.ir}
 \\
 {\small\em  School of Mathematics, Institute for Research in Fundamental Sciences (IPM)},\\{\small\em P.O. Box: 19395--5746, Tehran, Iran}}

\date{}

\begin{document}

\maketitle

%\fontsize{22}{23}\selectfont
%\baselineskip=12mm
\begin{abstract}
For a finite group $G$, let $\Delta(G)$ denote the character graph built on the set of degrees of the irreducible complex characters of $G$. In graph theory, a perfect graph is a graph $\Gamma$ in which the chromatic number of every induced subgraph $\Delta$ of $\Gamma$ equals the clique number of $\Delta$. In this paper,
 we show that the character graph $\Delta(G)$ of a finite group $G$ is always a perfect graph. We also prove that the chromatic number of the complement of $\Delta(G)$ is at most three.
 \end{abstract}
\noindent {\bf{Keywords:}}  Character graph, Character degree, Perfect graph. \\
\noindent {\bf AMS Subject Classification Number:}  20C15, 05C17, 05C25.

\section{Introduction}
$\noindent$ Let $G$ be a finite group . Also let ${\rm cd}(G)$ be the set of all character degrees of $G$, that is,
 ${\rm cd}(G)=\{\chi(1)|\;\chi \in {\rm Irr}(G)\} $, where ${\rm Irr}(G)$ is the set of all complex irreducible characters of $G$. The set of prime divisors of character degrees of $G$ is denoted by $\rho(G)$. It is well known that the
 character degree set ${\rm cd}(G)$ may be used to provide information on the structure of the group $G$. For example, Ito-Michler's Theorem \cite{[C]} states that if a prime $p$ divides no character degree of a finite group $G$, then $G$ has a
 normal abelian Sylow $p$-subgroup. Another result due to Thompson \cite{[D]} says that if a prime $p$ divides
 every non-linear character degree of a group $G$, then $G$ has a normal $p$-complement.

A useful way to study the character degree set of a finite group $G$ is to associate a graph to ${\rm cd}(G)$.
One of these graphs is the character graph $\Delta(G)$ of $G$ \cite{[I]}. Its vertex set is $\rho(G)$ and two vertices $p$ and $q$ are joined by an edge if the product $pq$ divides some character degree of $G$. We refer the readers to a survey by Lewis \cite{[M]} for results concerning this graph and related topics.

Let $\Gamma=(V(\Gamma),E(\Gamma))$ be a finite simple graph with the vertex set $V(\Gamma)$ and edge set $E(\Gamma)$.  A clique of $ \Gamma $ is a set of mutually adjacent vertices, and that the maximum size of a clique of $ \Gamma $, the clique number of $ \Gamma $, is denoted by $ \omega(\Gamma) $. Minimum number of colors needed to color vertices of the graph $ \Gamma $ so that any two adjacent vertices of $ \Gamma $ have different colors, is called the chromatic number of $ \Gamma $ and denoted by $ \chi(\Gamma) $. Clearly $\omega(\Gamma)\leqslant \chi(\Gamma)$ for any graph $\Gamma$. The graph $\Gamma$ is perfect if $\omega(\Delta)=\chi(\Delta)$, for every induced subgraph $\Delta$ of $\Gamma$.

The theory of perfect graphs relates the concept of graph colorings to the concept of cliques. Aside from having an interesting structure, perfect graphs are considered important for three reasons. First, several common classes of graphs are known to always be perfect. For instance, bipartite graphs, chordal graphs and comparability graphs are perfect. Second, a number of important algorithms only work on perfect graphs. Finally, perfect graphs can be used in a wide variety  of applications, ranging from scheduling to order theory to communication theory.

One of the turning point for the investigation on the character degree graph is the "Three-Vertex Theorem" by Palfy \cite{palfy}: the complement of $\Delta(G)$ does not contain any triangle whenever $G$ is a finite solvable group. In a recent paper  \cite{1} this was extended by showing that, under the same solvability assumption, the complement of $\Delta(G)$ dose not contain any cycle of odd length, which is equivalent to say that the complement of $\Delta(G)$ is a bipartite graph. Thus as the complement of a perfect graph is perfect \cite{lovasz}, the character graph $\Delta(G)$ of a solvable group $G$ is perfect. This argument motivates an interesting problem: is the character graph of an arbitrary finite group perfect?. In this paper, we wish to solve this problem.\\

\noindent \textbf{Theorem A.}  \textit{Let $G$ be a finite group then $\Delta(G)$ is a perfect graph.
}\\

Theorem A will be used to determine the chromatic number of the complement of $\Delta(G)$.\\

\noindent \textbf{Corollary B.}  \textit{ Suppose $G$ is a finite group. Then the chromatic number of the complement of $\Delta(G)$ is at most three.
}
%%%%%%%%%%%%%%%%%%%%%%%%%%%%%%%%%%%%%%
\section{Preliminaries}
$\noindent$ In this paper, all groups are assumed to be finite and all
graphs are simple and finite. For a finite group $G$, the set of prime divisors of $|G|$ is denoted by $\pi(G)$.  Also note that for an integer $n\geqslant 1$,  the set of prime divisors of $n$ is denoted by $\pi(n)$. We begin with Corollary 11.29 of \cite{[isa]}.

\begin{lemma}\label{fraction}
Let $ N \lhd G$ and $\chi \in \rm{Irr}(G)$. Let $\varphi \in \rm{Irr}(N)$ be a constituent of $\chi_N$. Then $\chi(1)/\varphi(1)$ divides $[G:N]$.
\end{lemma}

Let $\Gamma$ be a graph with vertex set $V(\Gamma)$ and edge set
$E(\Gamma)$. A cycle on $n$ vertices $v_1, \ldots,
v_n$, $n \geq 3$, is a graph whose vertices can be arranged in a
cyclic sequence in such a way that two vertices are adjacent if
they are consecutive in the sequence, and are non-adjacent
otherwise. A cycle with $n$ vertices is said to be of length $n$
and is denoted by $C_n$, i.e., $C_n : v_1, \ldots, v_n, v_1$.
 The join $\Gamma \ast \Delta$ of graphs $\Gamma$
and $\Delta$ is the graph $\Gamma \cup \Delta$ together with all
edges joining $V (\Gamma)$ and $V (\Delta)$.  An independent set is a set of vertices of a graph
$\Gamma$ such that no two of them are adjacent. A maximum
independent set is an independent set of largest possible size.
This size is called the independence number of $\Gamma$ and is
denoted by $\alpha(\Gamma)$. Finally we should mention that  the complement of $\Gamma$ and the induced subgraph of $\Gamma$ on $X\subseteq V(\Gamma)$
 are denoted by $\Gamma^c$ and  $\Gamma[X]$, respectively.  For more details, we
refer the reader to basic textbooks on the subject, for instance
\cite{BM}. Now we present some properties of perfect graphs.

\begin{lemma}\label{comp}\cite{lovasz}
A graph $\Gamma$ is perfect if and only if the complement of $\Gamma$ is perfect.
\end{lemma}

\begin{lemma}\label{hole}\cite{[cr]}
A graph $\Gamma$ is perfect  if and only if it has no induced subgraph isomorphic either to a cycle of odd order at least 5, or to the complement of such a cycle.
\end{lemma}

We now state some relevant results on character graphs
needed in the next section.

\begin{lemma}\label{pal}
Let $G$ be a group and let $\pi \subseteq \rho(G)$.\\
\textbf{a)} (Palfy's condition) \cite {palfy} If $G$ is solvable and $|\pi|\geqslant 3$, then there exist two distinct primes $u, v$ in $\pi$ and $\chi \in \rm{Irr}(G)$ such that $uv | \chi(1)$.\\
\textbf{b)} (Moreto-Tiep's condition) \cite{moreto} If $|\pi|\geqslant 4$, then there exists $\chi \in \rm{Irr}(G)$ such that $\chi (1)$ is divisible by two distinct primes in $\pi$.
\end{lemma}

\begin{lemma}\label{chpsl}\cite{[white]}
Let $G\cong \rm{PSL}_2(q)$, where $q\geqslant 4$ is a power of a prime $p$.\\
\textbf{a)}
 If $q$ is even, then $\Delta(G)$ has three connected components, $\{2\}$, $\pi(q-1)$ and $\pi(q+1)$, and each component is a complete graph.\\
\textbf{b)}
 If $q>5$ is odd, then $\Delta(G)$ has two connected components, $\{p\}$ and $\pi((q-1)(q+1))$.\\
i)
 The connected component $\pi((q-1)(q+1))$ is a complete graph if and only if $q-1$ or $q+1$ is a power of $2$.\\
ii)
  If neither of $q-1$ or $q+1$  is a power of $2$, then $\pi((q-1)(q+1))$ can be partitioned as $\{2\}\cup M \cup P$, where $M=\pi (q-1)-\{2\}$ and $P=\pi(q+1)-\{2\}$ are both non-empty sets. The subgraph of $\Delta(G)$ corresponding to each of the subsets $M$,$P$ is complete, all primes are adjacent to $2$, and no prime in $M$ is adjacent to any prime in $P$.
 \end{lemma}

\begin{lemma}\label{cycle} \cite{AC}
Let $G$ be a finite group, and let $\pi$ be a subset of the vertex set of $\Delta(G)$ such that $|\pi|$ is an odd number larger than 1. Then $\pi$ is the set of vertices of a cycle in $\Delta(G)^c$ if and only if $O^{\pi^\prime}(G)=S\times A$, where $A$ is abelian, $S\cong \rm{SL}_2(u^\alpha)$ or $S\cong \rm{PSL}_2(u^\alpha)$ for a prime $u\in \pi$ and a positive integer $\alpha$, and the primes in $\pi - \{u\}$ are alternately odd divisors of $u^\alpha+1$ and  $u^\alpha-1$.
\end{lemma}

 \begin{lemma}{\label{HIHLEW}}\cite{HIHL}
If $G$ and $H$ are two non-abelian groups that satisfy $\rho (G)
~\cap~ \rho (H) = F$ with $|F| = n$, then $\Delta (G\times H) = K_n \ast \Delta (G) [\rho (G) - F] \ast \Delta (H) [\rho (H) - F]$, where $K_n$ is a complete graph with vertex set $F$.
\end{lemma}
%%%%%%%%%%%%%%%%%%%%%%%%%
\section{Proof of main results}
\noindent In this section, we wish to prove our main results.\\

\noindent\textit{\bf Proof of Theorem A.}
On the contrary, we assume that $\Delta(G)$ is not perfect. Then by Lemma \ref{hole}, there exists a subset $\pi\subseteq \rho (G)$ such that for some integer $n\geqslant 2$, $|\pi|=2n+1$ and $\Delta(G)^c[\pi]$ or $\Delta(G)[\pi]$ is a cycle. Now one of the following cases occurs:\\
Case 1. $\Delta(G)^c[\pi]$ is a cycle. Then there exist primes $p_0, p_1, \dots , p_{2n} \in \rho(G)$ such that $\Delta(G)^c[\pi]$ is the cycle $C:p_0, p_1, \dots , p_{2n},p_0$. Using Lemma \ref{cycle},
 $N:=O^{\pi^\prime}(G)=S\times A$, where $A$ is abelian, $S\cong \rm{SL}_2(p^m)$ or $S\cong \rm{PSL}_2(p^m)$ for a prime $p\in \pi$ and a positive integer $m$, and the primes in $\pi - \{p\}$ are alternately odd divisors of $p^m+1$ and  $p^m-1$. Without loss of generality, we can assume that $p=p_0$. Since $\Delta(G)^c[\pi]$ is a cycle of length at least $5$, $p_1$ is not adjacent to both vertices $p_3$ and $p_4$ in $\Delta(G)^c$. Also for some $\epsilon \in \{\pm 1\}$, $p_3 \in \pi(p^m-\epsilon)-\{2\}$ and $p_4 \in \pi(p^m+\epsilon)-\{2\}$. Note that $p_1$ is an element of either $\pi(p^m-\epsilon)-\{2\}$ or $\pi(p^m+\epsilon)-\{2\}$.  Without loss of generality, we can assume that $p_1 \in\pi(p^m-\epsilon)-\{2\}$. since $p_1$ and $p_4$ are adjacent vertices in $\Delta(G)$, for some $\chi \in \rm{Irr}(G)$, $p_1p_4|\chi(1)$. Now let $\varphi \in \rm{Irr}(N)$ be a constituent of $\chi_N$. Then by Lemma \ref{fraction}, $\chi(1)/\varphi(1)$ divides $[G:N]$. Therefore as $G/N$ is a $\pi^\prime$-group, $p_1p_4|\varphi(1)$. It is a contradiction as using Lemma \ref{chpsl}, $p_1$ and $p_4$ are non-adjacent vertices in $\Delta(N)$.\\
Case 2. $\Delta(G)[\pi]$ is a cycle. If $\Delta(G)[\pi]\cong C_5$, then $\Delta(G)^c[\pi]\cong C_5$ which is a contradiction with Case 1. Thus $|\pi|\geqslant 7$. Hence there exist distinct primes $p_1, p_2, p_3,q_1, q_2, q_3\in \pi$ such that the induced subgraphs of $\Delta(G)^c$ on the sets  $\pi_1:=\{p_1, p_2, p_3\}$ and $\pi_2:=\{q_1, q_2, q_3\}$ are cycles of length $3$. Hence by Lemma \ref{cycle}, for every $i=1,2$, $N_i:=O^{\pi_i^\prime}(G)=S_i\times A_i$, where $A_i$ is abelian, $S_i\cong \rm{SL}_2(u_i^{\alpha_i})$ or $S_i\cong \rm{PSL}_2(u_i^{\alpha_i})$ for a prime $u_i\in \pi_i$ and a positive integer $\alpha_i$, and the primes in $\pi_i - \{u_i\}$ are alternately odd divisors of $u_i^{\alpha_i}+1$ and  $u_i^{\alpha_i}-1$. Let $N:=S_1S_2$. It is easy to see that $N/Z(N)\cong \rm{PSL}_2(u_1^{\alpha_1})\times\rm{PSL}_2(u_2^{\alpha_2})$. Therefore by Lemmas \ref{chpsl} and \ref{HIHLEW}, $\Delta(N/Z(N))[\pi_1 \cup \pi_2]\subseteq \Delta(G)$ is a complete bipartite graph with parts $\pi_1$ and $\pi_2$ which is a contradiction as $\Delta(G)[\pi]$ is a cycle of $\Delta(G)$.\qed\\

\noindent\textit{\bf Proof of Corollary B.}
Using Theorem A, $\Delta(G)$ is a perfect graph. Thus by Lemma \ref{comp}, $\Delta(G)^c$ is too. Hence $\chi(\Delta(G)^c)=\omega(\Delta(G)^c)$. It is clear that $\omega(\Delta(G)^c)=\alpha(\Delta(G))$. By Lemma \ref{pal}, $\alpha(\Delta(G))\leqslant 3$. Hence $\chi(\Delta(G)^c)\leqslant 3$ and the proof is completed.\qed
%%%%%%%%%%%%%%%%%%%%%%%%%%%%%%%%%%%%%%%%%%%%%%

\section*{Acknowledgements}
This research was supported in part
by a grant  from School of Mathematics, Institute for Research in Fundamental Sciences (IPM).

%%%%%%%%%%%%%%%%%%%%%%%%%%%%%

\end{document}